\documentclass[aip,
 jmp,
 amsmath  ,amssymb,
 reprint,%
]{revtex4-1}

\newcommand{\vek}[1]{\mathbf{#1}} %\underline{
\newcommand{\unitvek}[1]{\mathbf{\hat{#1}}}
\newcommand{\nn}[0]{\nonumber}

\newcommand{\mat}[1]{\underline{\underline{\mathbf{#1}}}}

\newcommand{\bfnabla}[0]{\boldsymbol{\nabla}}

\usepackage{graphicx}% Include figure files
\usepackage{dcolumn}% Align table columns on decimal point
\usepackage{amsbsy}
\usepackage{multirow}
\usepackage[abs]{overpic}
\usepackage{natbib}
\begin{document}

%\section*{Title: Solving Stress and Compliance Constrained Volume Minimization using Anisotropic Mesh Adaptation, the Method of Moving Asymptotes and a Global p-norm}
%\preprint{AIP/123-QED}
\title[Solving Stress Constrained Compliance Minimization using Anisotropic Mesh Adaptation, MMA and a p-norm]{Solving Stress Constrained Compliance Minimization using Anisotropic Mesh Adaptation, the Method of Moving Asymptotes and a Global p-norm}\author{Kristian Ejlebjerg Jensen}
\address{Imperial College London, Department of Earth Science and Engineering, London SW7 2AZ, United Kingdom}
\email{kristianejlebjerg@gmail.com}
\date{\today}

\begin{abstract}
\section*{abstract}
The p-norm often used in stress constrained topology optimisation supposedly mimics a delta function and it is thus characterised by a small length scale and ideally one would also prefer to have the solid-void transition occur over a small length scale, since the material in this transition does not have a clear physical interpretation. We propose to resolve these small length scales using anisotropic mesh adaptation. We use the method of moving asymptotes with interpolation of sensitivities, asymptotes and design variables between iterations. We demonstrate this combination for the portal and L-bracket problems with p=10, and we are able to investigate mesh dependence. Finally, we suggest relaxing the L-bracket problem statement by introducing a rounded corner.
\end{abstract}
%\pacs{Valid PACS appear here}
\keywords{Anisotropic, mesh, adaptation, topology, optimisation, Stress, Constraints, FEniCS, PRAgMaTIc.}
\maketitle

\section{Introduction}
Anisotropic mesh adaptation is an established technique for ensuring computational efficiency in the context of multiscale problems \cite{habashi2000anisotropic,pain2001tetrahedral,loseille2011continuous}, sometimes reducing the computational work with several orders of magnitude \cite{loseille2010fully}. Within the field of structural optimisation fixed structured meshes remain popular due to ease of implementation and compatibility with mathematical optimisers. It is thus mainly in the context of methods utilising various continuous sensitivities \cite{wallin2012optimal,amstutz2010topological} that mesh adaptivity has been applied, and even then mainly in the context of mesh refinement, but some of the most recent advancements also relate to the use of swapping operations\cite{christiansen2015combined}. Parallelism is another popular technique for speeding up computations \cite{borrvall2001large, aagetopology, aage2008topology} and for truly large scale problems, it indeed is the only option. 

The community of structural optimisation has devoted significant effort towards not only finding stiff and light structures, but also structures that do not break under load. That is the structure has to satisfy a stress constraint. Recent strategies include methods with global stress constraints, either with an explicit design representation, such that computation of void elements and related stress is avoided \cite{xia2012level}, or an implicit representation \cite{guo2011stress,zhang2013optimal}. An essential question in this context is a global versus a local constraint as investigated in the context of level-set methods \cite{guo2011stress,zhang2013optimal}. One can also choose to let the number of constraints vary throughout the optimisation \cite{bruggi2012topology}.

An older, but still popular approach \cite{duysinx1998new,holmberg2013stress}, approximates the local stress constraints with a single global constraint using a p-norm. The approximation is good for large values of p, but this causes numerical problems and therefore a compromise has to be made in practice. If the numerical problems are due to discretisation errors, one can hope to achieve a better compromise with more computational power or by being more efficient with the power available. It is the second option, that we choose to pursue with anisotropic mesh adaptation in this paper. The mesh adaptation introduces inconsistencies, but this drawback can be outweighed by the benefits of improved resolution of small length scales.

%The accuracy of continuous sensitivities is limited by mesh resolution, but so is the accuracy of the forward problem in the case of discrete sensitivities, so a mesh that minimises discretisation errors is desired regardless of the method. 
We choose to calculate discrete sensitivities as this allows us to harness the power of libraries for automatic adjoint derivation \cite{farrell2013automated}. A continuous sensitivity is required for driving the mesh adaptation as well as when interpolating between meshes,  but this is trivial to compute by dividing the discrete sensitivity with the design variable volumes, and one can invert the operation, if discrete variables are desired for the optimiser \footnote{Preliminary results related to this work was presented at FEniCS 14 in Paris, June 2014 and at the International Meshing Roundtable 23 in London, October 2014.}.

\section{Anisotropic Mesh Adaptation} \label{sec:adapt}
Stretched elements with small or large angles are normally discouraged in the context of the finite element method, as they give rise to not only large discretisation errors, but also cause problems with iterative solvers. This wisdom is, however, only applicable to problems without any anisotropic features in the solution. In fact, problems with strong anisotropy require anisotropic meshes for optimal use of the computational degrees of freedom. Note that the problems with large angles and iterative solvers are only severe for extreme anisotropy and, furthermore, it is possible to create anisotropic meshes with few large angles \cite{loseille2014metric}.

\begin{figure}[!htb]
\centering
\includegraphics[width=0.35\textwidth]{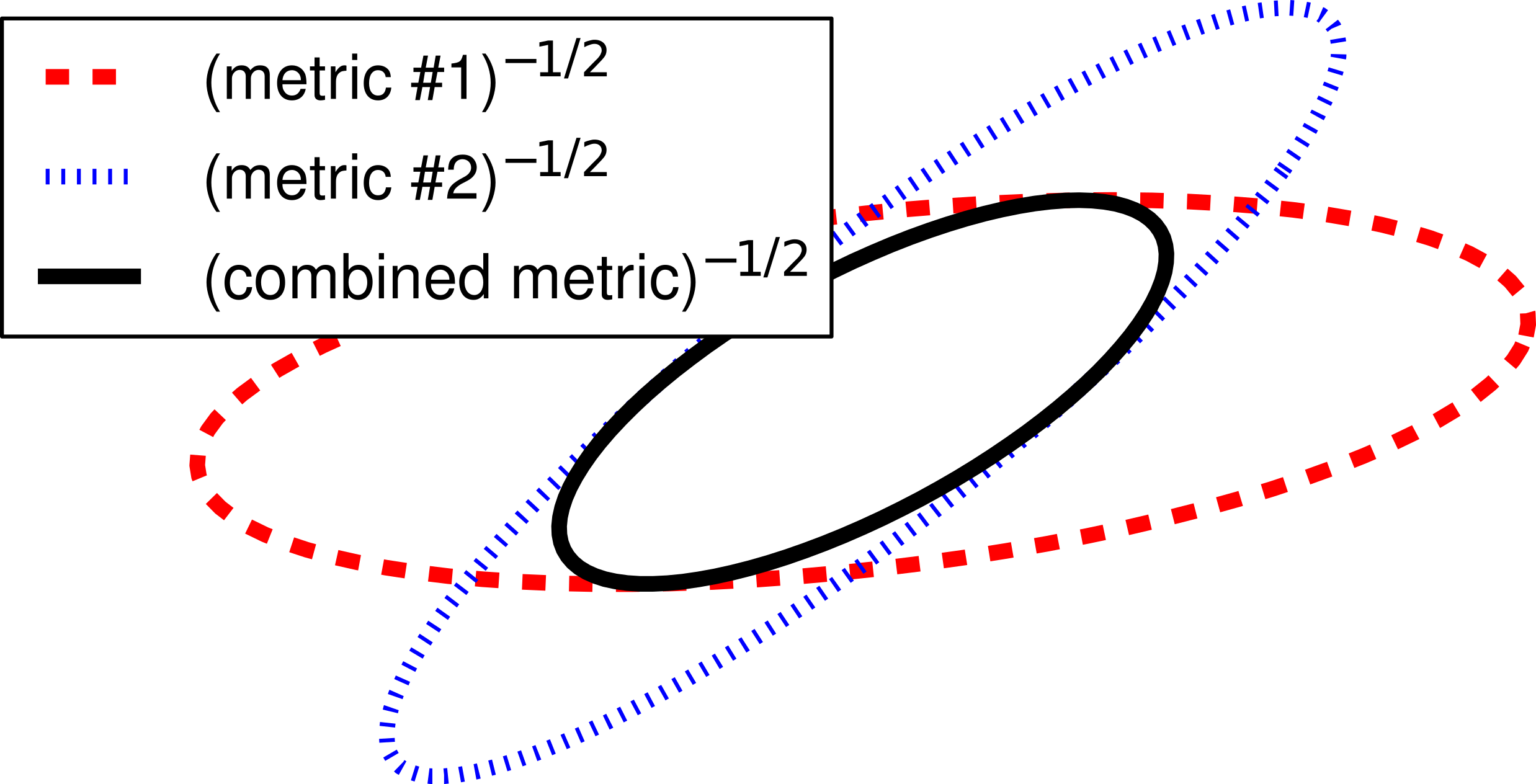} %,trim=0 0 0 5cm
\caption{The inner ellipse method is illustrated in the case of intersection, but it is common to see one ellipse entirely within the other, such that anisotropy is preserved.}\label{fig:ellipse}
\end{figure}

The properties of the quasi optimal mesh can be derived in a continuous sense \cite{loseille2011continuous} using a metric tensor, $\mat{\mathcal{M}}$. This is a symmetric and positive definite tensor field that maps the optimal element to the unit element, which has unit length edges. The first step is to compute the Hessian, $\mat{H}$, of the variable whose discretisation error is to be minimised and convert it into a positive definite matrix. This is done by taking the absolute value in the principal frame using the operator $\mat{\mathrm{abs}}$, which corresponds to removing the sign of the eigenvalues. The optimal mesh metric \cite{chen2007optimal} can then be expressed as 
\begin{eqnarray} 
\mat{\mathcal{M}} = \frac{1}{\eta}\left(\mathrm{det}[\mat{\mathrm{abs}}(\mat{H})]\right)^{-\frac{1}{2q+d}} \mat{\mathrm{abs}}(\mat{H}) \label{eqn:M} ,  
\end{eqnarray}
where $d$ is the number of dimensions, $q$ is the error norm to be minimised, $\mathrm{det}$ is the determinant and $\eta$ is a scaling factor. It is possible to combine the metrics of several variables using the inner ellipse method illustrated in figure \ref{fig:ellipse}, see \cite{pain2001tetrahedral} for implementation details. Note that the unit of $\eta$ is always so that the metric in equation (\ref{eqn:M}) has dimension of squared inverse length.

\begin{figure}[!htb]
\centering
\vspace{1mm}
\includegraphics[width=0.45\textwidth]{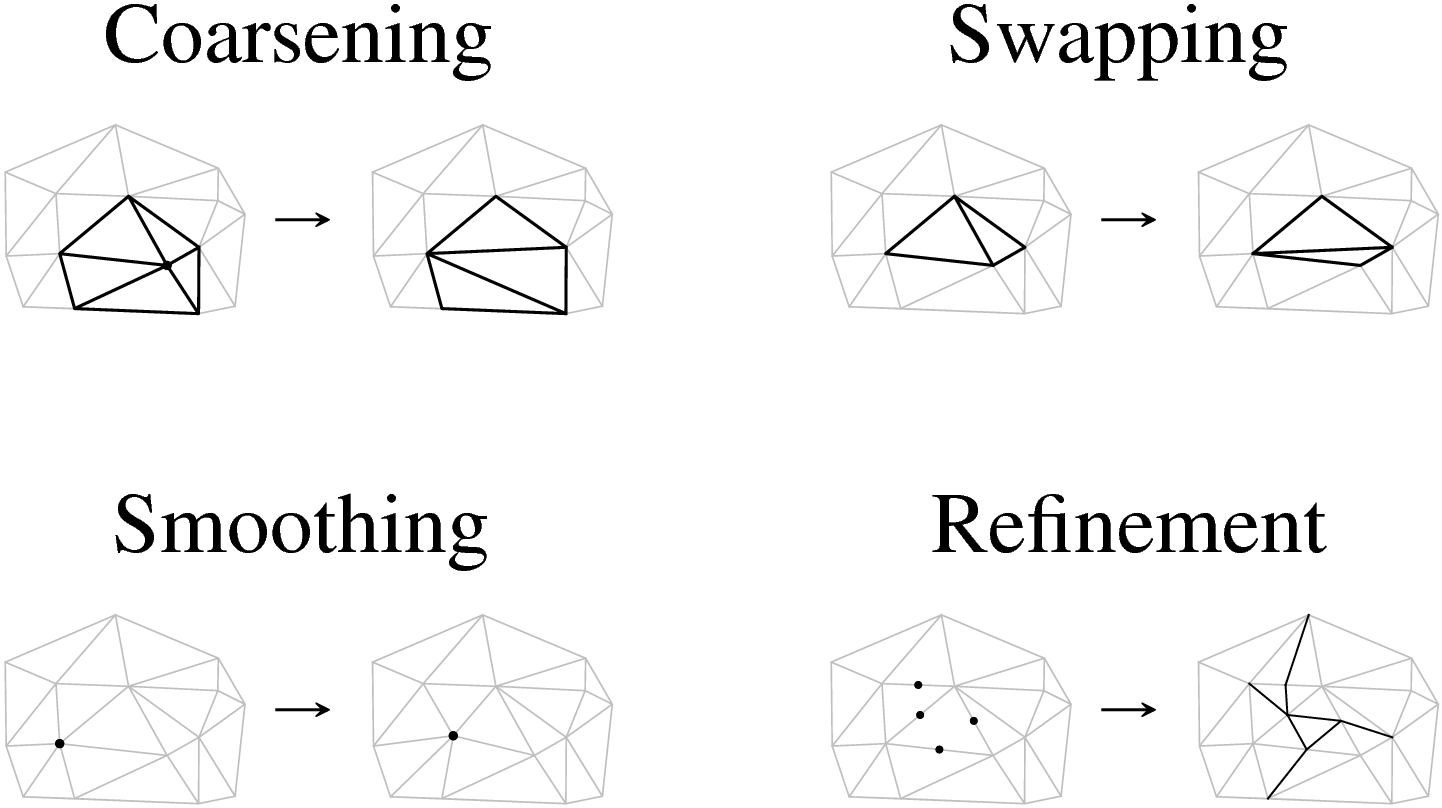}
\caption{Four local mesh modifications are illustrated: Coarsening, refinement, swapping and smoothing. Only the refinement operation is allowed to increase the worst local element quality for the MATLAB/Octave implementation, while the PRAgMaTIc mesh generator also allows the worst element quality to deteriorate for the coarsening.}\label{fig:meshmod}
\end{figure}

Once a final metric has been calculated, it can be passed to an anisotropic mesh generator together with the current mesh, in fact our metric is defined on the nodes of the current mesh. We use two mesh generators, which applies four local mesh modifications: Coarsening, refinement, swapping and smoothing, see figure \ref{fig:meshmod}. The mesh generators work by applying the modifications to optimise a heuristic quality measure\cite{vasilevski2005error}, which quantifies the difference between the discrete mesh and the optimal continuous one. We use an optimised C++ mesh generator, PRAgMaTIc \cite{rokos2013thread} to illustrate the negligible performance penalty of mesh adaptation, but presently this does not support curved geometries, so in that case we resort to a slower MATLAB/Octave implementation.

Our implementation can generate a metric from a function represented on continuous 2nd order polynomials, but 1st order polynomials can also be used. In the latter case Galerkin projection is used as derivative recovery technique to compute the Hessian. In any case, Galerkin projection is used to convert the element wise constant metric to a node based representation, before it is passed to the mesh generator.

\section{Topology Optimisation}
In an effort to maximise the impact of this paper, we have tried to choose the most popular method for topology optimisation with stress constraints, but we see no reason why other methods could not benefit equally well from the use of mesh adaptation. Similarly other problems with small length scales, such as many convective problems, could most definitely be solved more efficiently using anisotropic mesh adaptation. This being said, we will focus on the density method with SIMP penalisation \cite{bendsoe2003topology} and a p-norm for relaxing the local stress constraints to a single global constraint \cite{duysinx1998new}. The idea is that the local constraint is satisfied in the limit of $p$ going to infinity. In practice a finite value is chosen, but in the following we will show that adaptive meshing allows for $p=10$ and that this effectively removes kinks in the design. We use a stress penalisation scheme to eliminate problems with void stress \cite{cheng1997varepsilon}, and the method of moving asymptotes is used to update the design variables \cite{svanberg1987method}.

\subsection{Mesh Dependence}
It is likely that the perfect solution to stress constrained topology optimisation will continue to allude the scientific community, but one might hope that scholars demonstrate their various heuristic schemes in the context of mesh %\footnote{It is not uncommon to see papers where convergence in the optimisation dimension is demonstrated in the context of the filter with a minimum length scale identical to the mesh resolution.} 
 independence such that the uncertainties of objective and constraint functions can be estimated, allowing for techniques to be compared on an objective basis. We have found the selection of papers including quantitative information related to computational cost in the context of stress constrained topology optimisation to be rather scarce\cite{amstutz2010topological,bruggi2012topology,suresh2013stress}, and we have been unable to identify any papers containing a discussion of mesh independence. This might be due to it being well known that the stress does not converge for designs with kinks, but in fact it has been shown that the design curvature can be effectively controlled using projection methods \cite{wang2011projection}.

It is well known that the p-norm underestimates the maximum stress, but we have deliberately chosen not to implement a ''local stress fix'' scheme \cite{le2010stress,suresh2013stress}, that is a scheme that changes the maximum stress to a conservative value such that the local stress constraint is satisfied. This is due to the fact that local stress constraints constitute a much more difficult problem. The alternative is to solve the relaxed stress problem first, and address the local stress problem by employing a conservative maximum stress in combination with a large p-norm. Regardless, the design has to go through a manual post-processing step as it is already typical for designs obtained with topology optimisation, but the topology is normally fixed in the post-processing and therefore it is important that the optimisation result is robust and thus mesh independent.

\section{Problem Setup} \label{sec:setup}

\begin{figure}[!htb]
\centering
\includegraphics[width=0.35\textwidth,clip,trim=2cm 13mm 0 0]{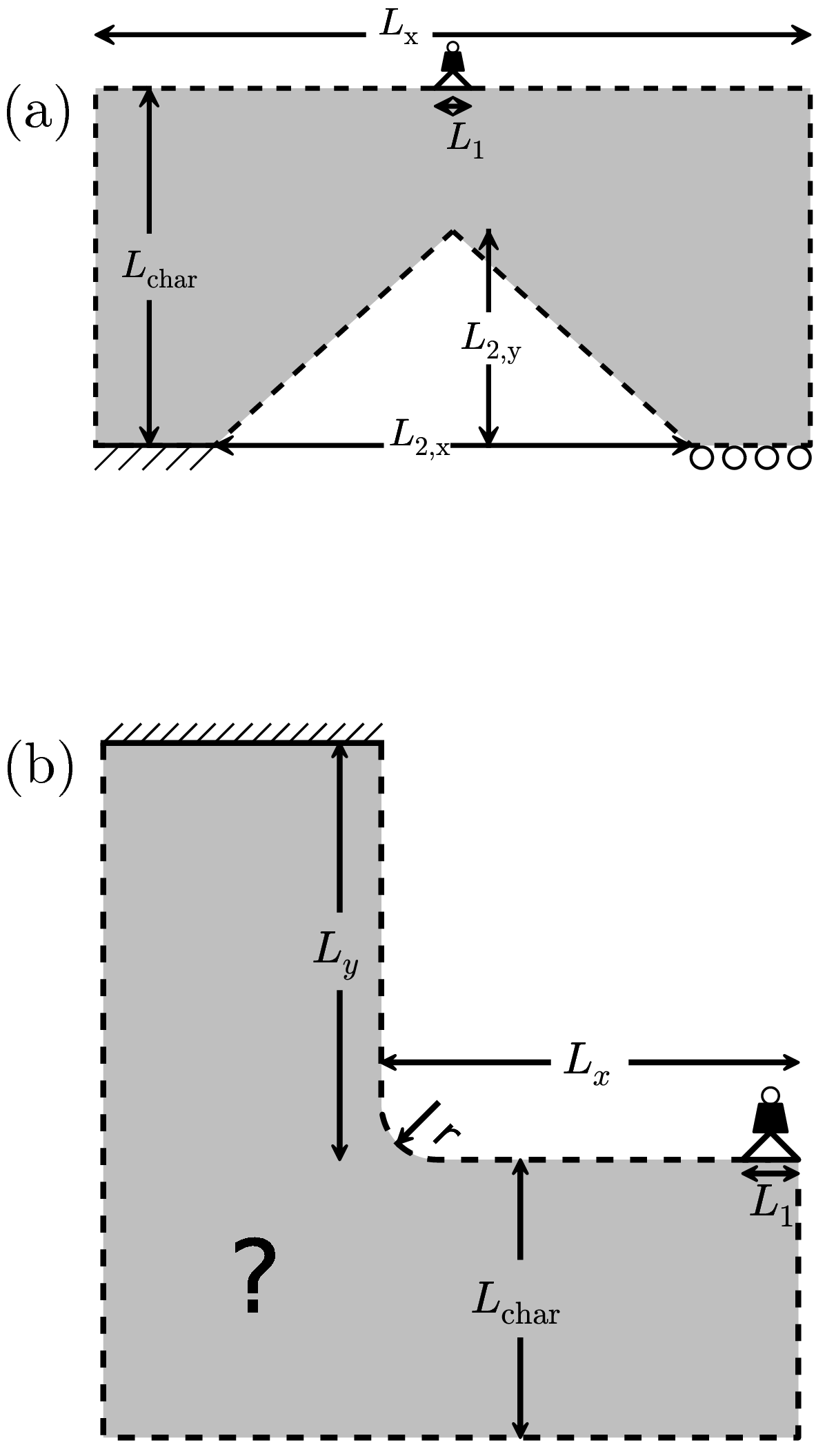}
\caption{Both the portal (a) and the L-bracket (b) are good benchmarks for stress constrained topology optimisation due to the stress concentration in the corners of the geometries. The load is distributed over the length $L_1$ in order to avoid a stress singularity in the problem definitions. Sharp corners can also give rise to stress singularities, and we have investigated the effect of this by rounding the corner of the L-bracket. The rollers to the lower right of the portal frame are associated with zero normal displacement, $\Omega_\mathrm{ice}$}\label{fig:Lgeom}
\end{figure}

We consider two 2D problems, the L-bracket and the portal frame \cite{le2010stress}, both with finite load and support areas, $\Omega_\mathrm{load}$ and $\Omega_{u=0}$, as illustrated in figure \ref{fig:Lgeom}. We model plane stress and use a Helmholtz filter \cite{lazarov2011filters} to compute a filtered design, $\tilde{\rho}$, with a minimum length scale, $L_\mathrm{min}$,
\begin{eqnarray}
\vek{0} &=& \bfnabla \cdot \mat{\sigma}, \quad \mat{\sigma} = 2G \mat{\epsilon}+\lambda \mat{I}\left(\mathrm{Tr} (\mat{\epsilon})+\partial_z u_z\right),  \label{eqn:stress} \\
\mat{\sigma} &=& \mat{\sigma}_\mathrm{load} \quad \mathrm{at} \quad \Omega_\mathrm{load}, \quad \vek{u}\cdot\vek{n} = 0 \quad \mathrm{at} \quad \Omega_\mathrm{ice} \nn \\
\vek{u}&=&\vek{0} \quad \mathrm{at} \quad \Omega_{u=0} \nn \\
\mat{\epsilon} &=& \frac{1}{2}\left(\bfnabla \vek{u} + [\bfnabla \vek{u}]^T \right), \quad \partial_z u_z = -\frac{\nu}{1-\nu} \bfnabla\cdot\vek{u}, \nn \\
G &=& E\frac{1}{2(1+\nu)}, \quad \lambda = E\frac{\nu}{(1+\nu)(1-2\nu)}, \nn \\
%\tilde{\rho} &=& \rho + L_\mathrm{min}^2\nabla^2\tilde{\rho} \label{eqn:filter}\\
%\tilde{\rho} &=& \rho + \bfnabla\cdot \mat{L}_\mathrm{min}^2 \cdot \bfnabla \tilde{\rho} \label{eqn:filter}\\
\tilde{\rho} &=& \rho + \bfnabla\cdot L_\mathrm{min}^2 \cdot \bfnabla \tilde{\rho} ,\label{eqn:filter}\\
E &=& E_\mathrm{min}+(E_\mathrm{max}-E_\mathrm{min}) \rho^{P_E} \label{eqn:E}
\end{eqnarray}
where $\vek{u}$ is the Lagrangian displacement, $\mat{\sigma}$ is the stress, $\mat{\epsilon}$ is the strain, $\partial_z u_z$ is the out of plane deformation, $\mat{I}$ is the identity tensor, $\mathrm{Tr}$ is the trace, $G$ is the shear modulus, $\lambda$ is Lam{\'e}'s  first parameter, $E$ is Young's modulus, $\nu$ is the Poisson ratio, $\rho$ is the design variable, and $P_E$ is the SIMP penalisation exponent. Note that we intend to use a sensitivity filter for the compliance. %, which corresponds to minimizing the compliance of a nonlocal elasticity problem \cite{sigmund2012sensitivity}. 
We thus avoid using the filtered design variable in equation (\ref{eqn:E}), it is only calculated for the purpose of driving the mesh adaptation.

The displacement as well as the design variables are discretised with continuous 1st order polynomials, while we use 2nd order when applying the PDE filter (\ref{eqn:filter}). The use of discontinuous constant design variables gives rise to excessive stress concentrations in the context of a sensitivity filter, and thus we have found significantly better results (objective functions) using nodal densities. Element wise constant densities might be used with a density filter, but then the issue of dealing with negative filtered design variables in a robust way arises, hence we use a sensitivity filter. 

The forward problem defined in equations (\ref{eqn:stress}-\ref{eqn:E}) is solved using FEniCS, an open source finite element engine with a high degree of automation \cite{LoggMardalEtAl2012a}. We use a direct solver for the forward, adjoint and filter problems, and an iterative solver for Galerkin projections between different element types.

\begin{figure}[!htb]
\centering
\includegraphics[width=0.4\textwidth]{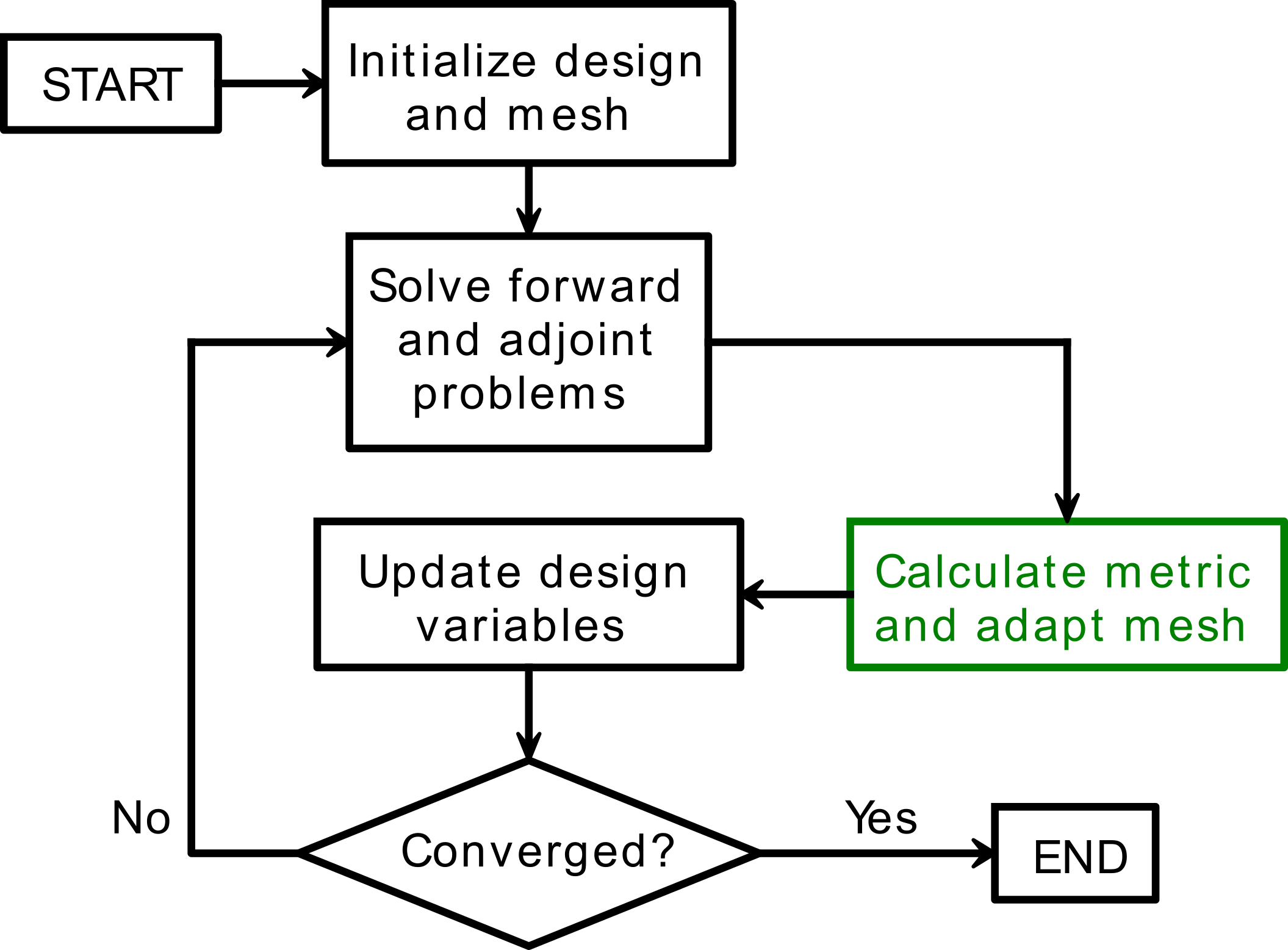}
\caption{This flowchart shows the position of the mesh adaptation in the design optimisation procedure, between the adjoint problem and the optimiser.}\label{fig:flowchart}
\end{figure}

We use the current design and displacement to calculate objective and constraint functions,
\begin{eqnarray}
%\sigma_\mathrm{miss} &=& E_S \sqrt{\frac{(\epsilon_{11}-\epsilon_{22})^2+\epsilon_{11}^2 +\epsilon_{22}^2+6\epsilon_{12}^2}{2}} \nn \\ %\label{eqn:s} \\
\sigma_\mathrm{miss} &=& E_S \sqrt{\epsilon_{xx}^2+\epsilon_{yy}^2-\epsilon_{xx}\epsilon_{yy}+3\epsilon_{12}^2} \nn \\ %\label{eqn:s} \\
E_S &=& E_\mathrm{min} + (E_\mathrm{max} - E_\mathrm{min}) \rho^{P_S} \label{eqn:ES}\\
V &=& \int_\Omega \rho d\Omega \nn \\ %\label{eqn:V} \\
C &=& \frac{1}{C_\mathrm{max}} \int_{\partial \Omega_\mathrm{load}} \vek{u} \cdot \mat{\sigma}_\mathrm{load} \cdot \unitvek{n} ds  - 1 \nn \\ %\label{eqn:C}\\
S &=& \left[\int_\Omega \left(\sigma_\mathrm{miss}/\sigma_\mathrm{max}\right)^p\right]^{1/p} - 1 \nn , %\label{eqn:S} ,
\end{eqnarray}
where $C_\mathrm{max}$ is the maximum compliance, $P_S$ is the stress penalisation exponent, $\sigma_\mathrm{miss}$ is the von mises stress and $\sigma_\mathrm{max}$ is its maximum value. The interpolations of the Young's modulus in equations (\ref{eqn:E}) and (\ref{eqn:ES}) agree in the limits of $\rho=0$ and $\rho=1$ corresponding to void and material, respectively, but one has to penalise the intermediate values differently for stress constrained optimisations \cite{cheng1997varepsilon}. The discrete gradient of $V$ can be calculated explicitly, and we use dolfint-adjoint \cite{farrell2013automated} to calculate the gradients of $S$ and $C$. These two discrete gradients are converted to continuous ones by division with the gradient of $V$ and anisotropic Helmholtz smoothing is applied to the stress sensitivity. This smoothing is defined on the continuous level using an equation identical to (\ref{eqn:filter}) with a tensor version of $L_\mathrm{min}$, which is based on the Steiner ellipse of the elements. We use the Helmholtz filter (\ref{eqn:filter}) for the compliance sensitivity,  %$C_\rho$,
%%
%\begin{eqnarray}
%%\tilde{C_\rho} &=& C_\rho\rho + \bfnabla\cdot L_\mathrm{min}^2 \cdot \bfnabla \tilde{C_\rho} ,\nn \\
%%\hat{C_\rho} &=&  \tilde{C_\rho}/\rho \nn
%\hat{C_\rho} &=& C_\rho\rho + \bfnabla\cdot L_\mathrm{min}^2 \cdot \bfnabla \hat{C_\rho} ,\nn 
%\end{eqnarray}
%%
%where $\hat{C_\rho}$ is the filtered compliance sensitivity, 
that is we do not use a sensitivity filter related to non-local elasticity as discussed in \cite{sigmund2012sensitivity}.
%The compliance sensitivity is multiplied by the design variables; then the minimum length scale is imposed using a Helmholtz filter (equation (\ref{eqn:filter)) and finally it is divided with the design variables.
 These continuous and smooth sensitivities are then used to calculate metrics associated with the compliance and stress constraints.  The filtered design variable is used to calculate a metric associated with the volume constraint and all three metrics are combined using the inner ellipse method, see figure \ref{fig:ellipse}. 
After the metric has been used to adapt the mesh, the optimiser variables (the asymptotes, continuous sensitivities, previous, and previous previous design variables) are interpolated on to the new mesh. Due to extrapolation at curved boundaries it can become necessary to enforce the box constraints on the asymptotes and design variables after this interpolation step.
Once this interpolation step is complete, the continuous sensitivities are converted to discrete ones, and the optimiser is called to update the design variables. This completes the optimisation loop as illustrated in figure \ref{fig:flowchart}. This methodology is inconsistent from a mathematical programming point of view, but except for the compliance sensitivity filter, we expect the errors to decrease with mesh refinement.

We consider volume minimisation under stress and compliance constraints, the formal problem statement being 
\begin{eqnarray}
\stackrel{\mathrm{min}}{_{0\leq\rho\leq 1}} V&,& \quad s.t. \nn \\
C \leq 0&,& \quad S \leq 0 \quad \mathrm{and} \quad \mathrm{eq.}\:\mathrm{ (\ref{eqn:stress}-\ref{eqn:E})} \nn 
\end{eqnarray}
%Note that we set a non-zero lower bound, $\rho_\mathrm{min}$ for the design variable in order to avoid negative filtered design variables due to the PDE filter in equation (\ref{eqn:filter}). 
The optimisation is initialised with a completely solid design, $\rho_0=1$. 

We make the problem dimensionless using $L_\mathrm{char}$ and $E_\mathrm{max}$ as characteristic length scale and stress, respectively. This is reflected in the set of parameters used in the optimisations,
\begin{eqnarray}
L_1 &=& 0.1 L_\mathrm{char}, \quad \mat{\sigma}_\mathrm{load} = E_\mathrm{max}/L_\mathrm{char}, \quad  \nu=0.3, \nn \\
L_\mathrm{min} &=& 5\cdot 10^{-2}L_\mathrm{char}, \quad \sigma_\mathrm{max} = 1.5 E_\mathrm{max},\quad p = 10, \nn \\
P_E &=& 3, \quad P_S = 0.5, \quad \mathrm{and} \quad E_\mathrm{min} = 10^{-3} E_\mathrm{max} . \nn
\end{eqnarray}
The maximum compliance and two of the length scales differ between the two benchmarks as shown in table \ref{tab:tab1}. Ideally the values for the maximum stress and compliance should be chosen according to the application, but for the academic benchmarks considered here, we just choose a maximum compliance that results in a non-extreme volume fraction and a maximum stress that results in the stress constraint being active.

\begin{table}[!htb]
\centering
\begin{tabular}{c | c c}
parameter \textbackslash \, geometry  & portal & L-bracket \\ \hline
$L_x/L_\mathrm{char}$ & 2 & 1.5 \\
$L_y/L_\mathrm{char}$ & 1 & 1.5 \\
$C_\mathrm{max}L_\mathrm{char}^{-2}/E_\mathrm{max}$ & 0.75 & 2.5 \\
$L_{2,x}/L_\mathrm{char}$ & 2-1/6 &  \\
$L_{2,y}/L_\mathrm{char}$ & 0.6 &  \\
\end{tabular}
\caption{The $L_x$ and $L_y$ length scale differs between the portal and L-bracket problems, and the maximum compliance is significantly higher for the L-bracket problem.}\label{tab:tab1}
\end{table} 

In addition to these we have the scaling factor, $\eta$, and four numerical parameters,
\begin{eqnarray}
c_\mathrm{MMA} &=& 10^3, \quad q=2, \quad h_\mathrm{min}=10^{-3}, \quad \mathrm{and} \nn \\
\mathrm{AR}_\mathrm{max} &=& 50 , \nn
\end{eqnarray}
where the $c_\mathrm{MMA}$ is related to enforcement of constraints, $q$ is the error norm to be minimised, $h_\mathrm{min}$ is the minimum element edge length and $\mathrm{AR}_\mathrm{max}$ is the maximum element aspect ratio. Finally, we use move limits, $\Delta \rho$, for the optimiser, which enforce the constraints
\begin{eqnarray}
\mathrm{abs}(\rho_{i+1}-\rho_{i}) &\leq& \Delta \rho , \nn
\end{eqnarray}
and we fix these at $\Delta \rho=0.1$.
%
%To avoid locking of the design, we widen the asymptotes, $\rho_\mathrm{low}$ and $\rho_\mathrm{upp}$,
%%
%\begin{eqnarray}
%\rho_\mathrm{low} & \equiv & \mathrm{max}\left(\rho_\mathrm{low}-4\Delta \rho(\rho-\rho_\mathrm{min})(1-\rho),\rho_\mathrm{min} \right) \nn \\
%\rho_\mathrm{upp} & \equiv & \mathrm{min}\left(\rho_\mathrm{upp}+4\Delta \rho(\rho-\rho_\mathrm{min})(1-\rho),1 \right) \nn
%\end{eqnarray}
%%
%This is turned off for the last 40 iterations, and the mesh adaptation is turned off for the last 20 iterations.

\section{Results}
%We fix the MMA $c$ parameter at $10^3$. and calculate mesh metrics corresponding to minimisation of the 2-norm ($q=2$ in equation (\ref{eqn:M})) for all mesh metrics. We impose a minimum edge length of $10^{-3}L_\mathrm{char}$ and a maximum aspect ratio of 50, but there is no constraint on the maximum number of elements.%for the elements, nor a maximum edge length or aspect ratio. % and set a minimum edge length of $10^{-6}$ in all optimisations.
In order to investigate mesh independence, we choose the scaling factor related to the metric of the compliance sensitivity as the primary numerical parameter to be varied and scale the number of iteration $\mathrm{it}_\mathrm{max}$, 	the move limits and the other scaling factors with a dimensionless version of this, $\eta_{\tilde{\rho}}$,
\begin{eqnarray}
\mathrm{it}_\mathrm{max} &=& \mathrm{\mathbf{round}}\left({600\sqrt{0.015/\eta_{\tilde{\rho}}}}\right) \nn \\
%\Delta \rho &=& 0.1\sqrt{\eta_{\tilde{\rho}}/0.015} \nn \\
\eta_C &=& \eta_{\tilde{\rho}} \nn \\
\eta_S &=& 4\eta_{\tilde{\rho}} , \label{eqn:etas}
\end{eqnarray}
%
%\paragraph{Compliance minimization}
where the $\mathrm{\mathbf{round}}$ function rounds the input to the nearest integer. Note the factor of four in equation (\ref{eqn:etas}), which reflects a lower emphasis on resolving the sensitivity of the stress constraint. It is possible to optimise with $\eta_S=\eta_{\tilde{\rho}}$, but this increases the tendency of the method to produce mesh dependent topologies. This observation and the factor of four is probably sensitive to the choice of $p=10$.

The sensitivity filter and interpolation of internal optimiser variables between meshes introduces inconsistencies in our implementation, and therefore we do not expect convergence in a strict sense, so we just plot the design variables and mesh elements for the iteration corresponding to the smallest volume at which the constraints are satisfied to the tolerance of the MMA $c$ parameter. The design topology in the final iteration is, however, identical to the one shown\footnote{except for the one case where an optimisation of the rounded L-bracket problem fails to recover from an infeasible design.}.

\begin{figure}[!htb]
\centering
\includegraphics[width=0.45\textwidth]{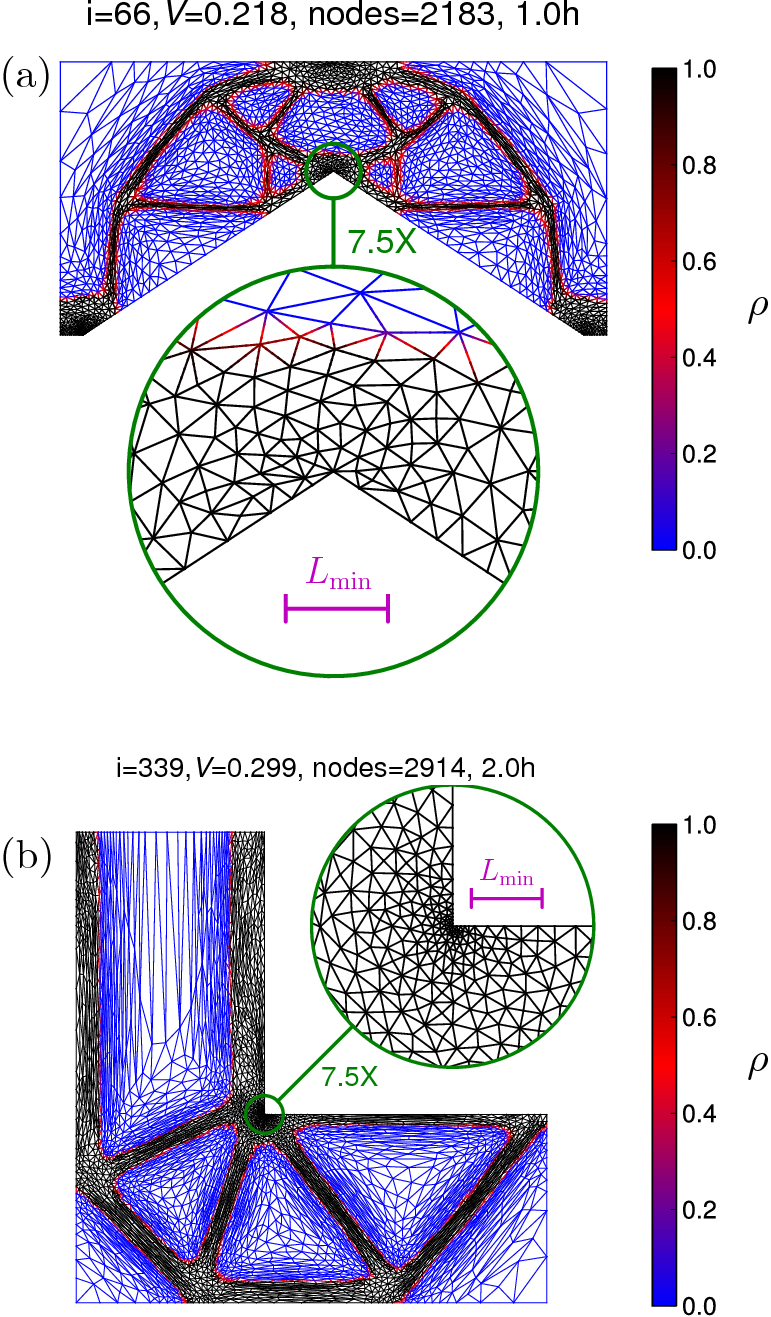}
\caption{The result of optimisations using $\sigma_\mathrm{max}=15E_\mathrm{max}$ and $\eta_{\tilde{\rho}}= 0.03$ are plotted. The large maximum stress results in designs similar to compliance minimisation, thus the bar going horisontally from the load in the L-bracket.}\label{fig:compliance}
\end{figure}

For reference we have performed optimisations with a large maximum stress ($\sigma_\mathrm{max}=15E_\mathrm{max}$) to mimic the results of pure compliance minimisation. This is shown in figure \ref{fig:compliance}, which shows the design kinks at the concave corners that we wish to eliminate by imposing the stress constraint.

\begin{figure*}[!htb]
\centering
\includegraphics[width=\textwidth]{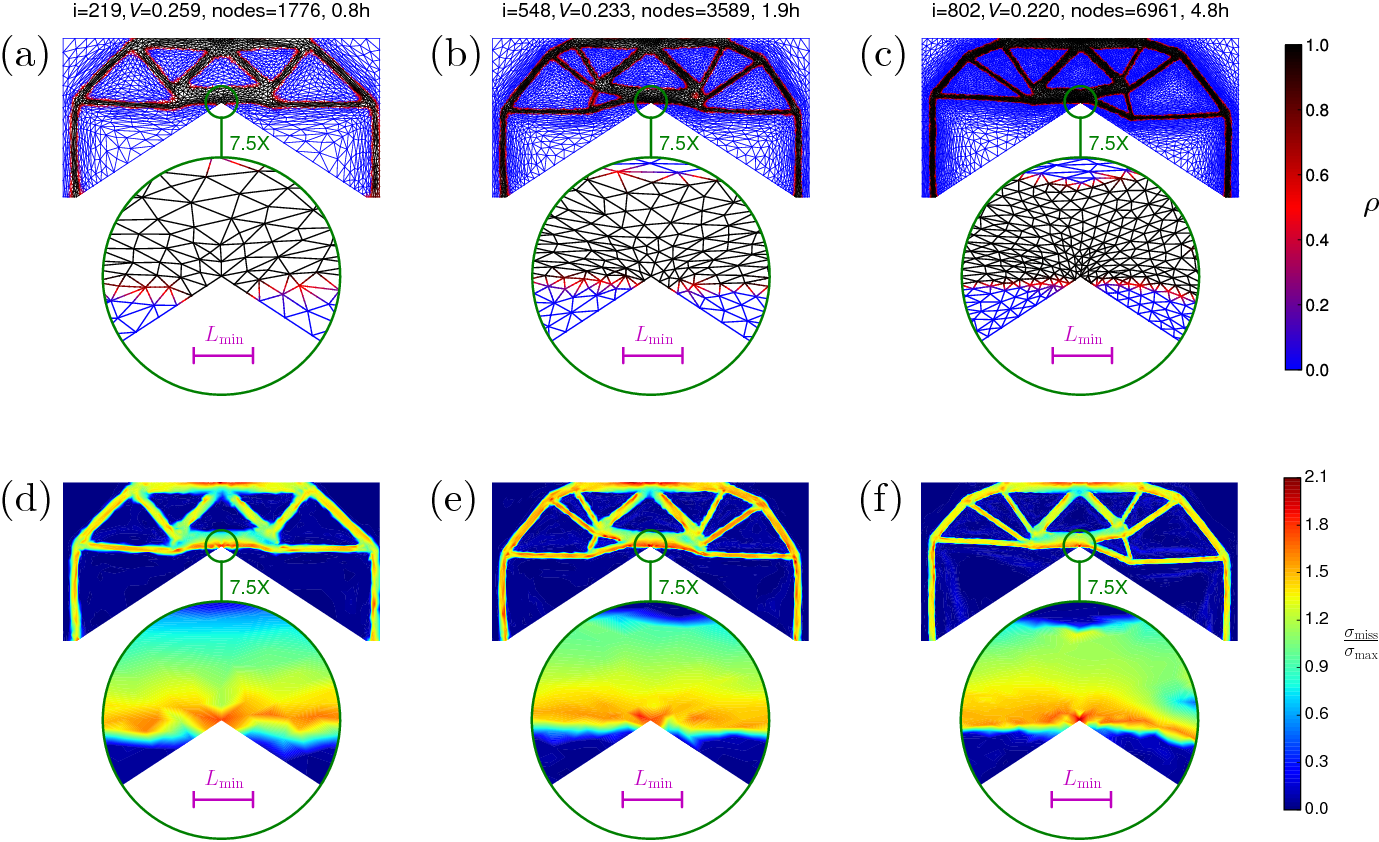}
\caption{The portal problem is optimised with $\eta_{\tilde{\rho}}=0.03$ (a,d), $0.015$ (b,e) and $0.0075$ (c,f). Note the local minimum for the finer optimisations, which is indicative of local minima.}\label{fig:portal}
\end{figure*}

For the portal problem, we have performed optimisations with $\eta_{\tilde{\rho}}=0.03$, $0.015$ and $0.0075$ as shown in figure \ref{fig:portal}. With continuous and linear displacements we get element wise constant strains, so we have to use Galerkin projection to compute continuous and linear stresses, as shown in \ref{fig:portal}(d-f). The finer optimisations most likely find local asymmetric minima, but in all cases the design kink is eliminated.

\begin{figure}[!htb]
\centering
\includegraphics[width=0.45\textwidth]{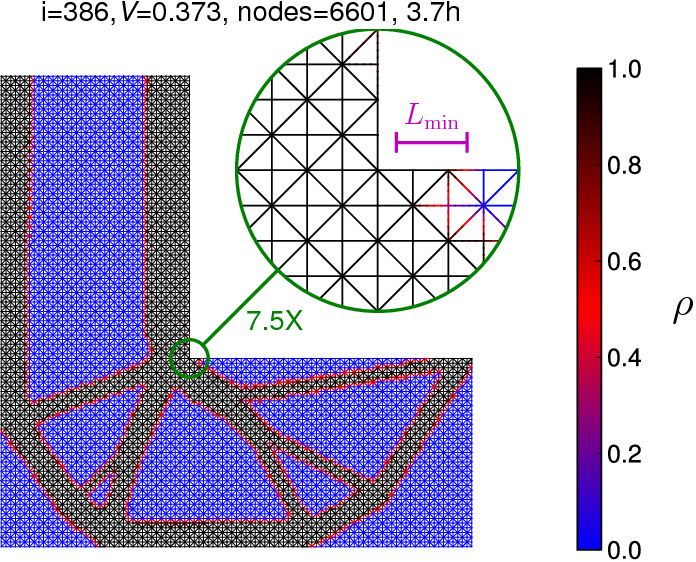}
\caption{The result of stress minimisation on a fixed mesh is shown for comparison with the results based on adaptive meshing.}\label{fig:static}
\end{figure}

\begin{figure*}[!htb]
\centering
\includegraphics[height=0.9\textheight]{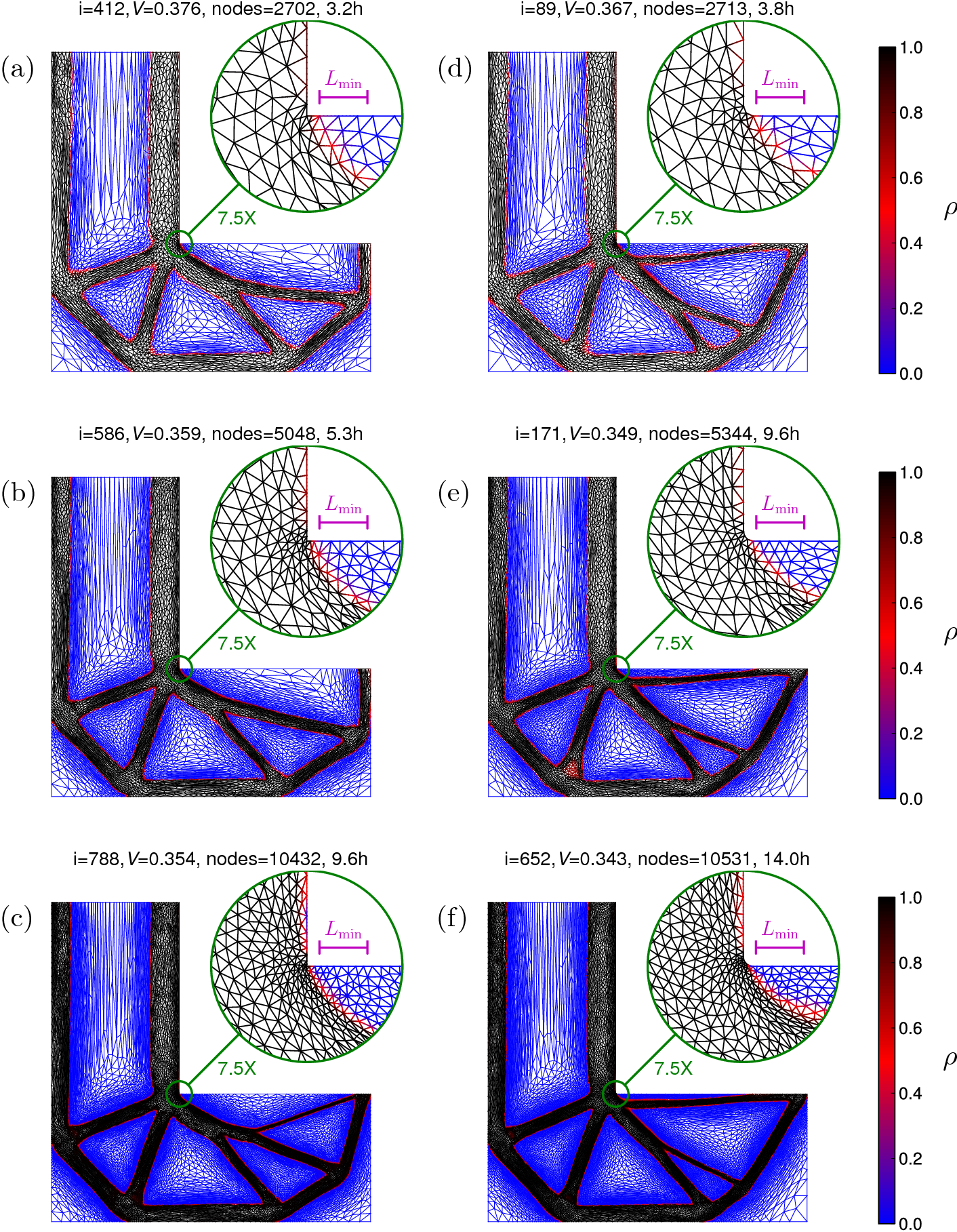} %3
\caption{Optimisations with a sharp (left) and a rounded corner (right) for $\eta_{\tilde{\rho}}$ equal to 0.03 (top), 0.015 (middle) and 0.0075 (bottom). The design variables and mesh elements are shown for the iteration (i) at which the lowest volume fraction ($V$) occurs, while the stress and compliance constraints are satisfied. See online version for colours.}\label{fig:designs}
\end{figure*}

For the L-bracket problem, we have performed an optimisation with a fixed mesh and $\mathrm{it}_\mathrm{max}=600$, as shown in figure \ref{fig:static}. This is to serve as a reference for the optimisations based on adaptive meshing. These are performed with a sharp ($r=0$) and a rounded corner ($r=0.01L_\mathrm{char}$) for $\eta_{\tilde{\rho}}$ equal to 0.03, 0.015 and 0.0075 as shown in figure \ref{fig:designs}. 
The two coarser optimisations with a sharp corner have a component in compression at the load, which might be unstable to perturbations in the load direction. This behaviour has also been observed previously\cite{holmberg2013stress,suresh2013stress}, and most likely it can be fixed by using a second load case. It takes more iterations to find the best feasible design for the sharp corner, possible indicating that the optimiser has an easier time dealing with the rounded corner, which also might explain the mesh dependent designs for the sharp corner. %Furthermore the topology of the optimisations with a rounded corner are not mesh independent, which could probably be addressed by some kind of continuation scheme.

\begin{figure}[!htb]
\centering
\includegraphics[width=0.45\textwidth]{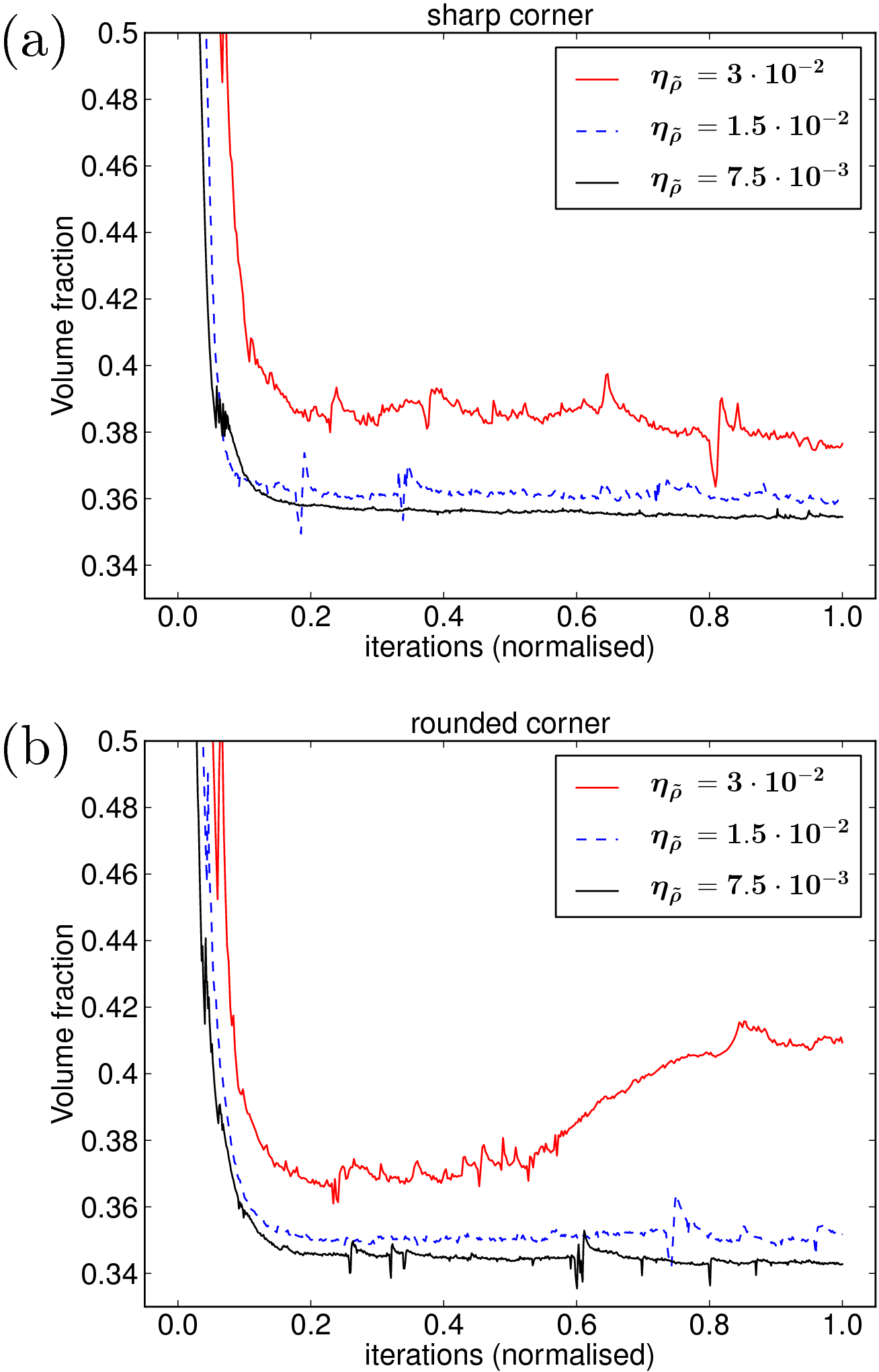}
\caption{The volume fraction is plotted versus normalised iteration numbers for different values of $\eta_{\tilde{\rho}}$ in the case of a sharp (a) as well as a rounded corner (b). In all cases there are some oscillations, which are probably caused by infeasible designs.}\label{fig:vol}
\end{figure}

It is well known that structures become weaker as the mesh is refined, i.e. the compliance converges from below \cite{lazarov2011filters} and the same is true for the stress. One would thus expect the volume to converge from below in a stress and compliance constrained optimisation problem, but this is not what we see in figure \ref{fig:vol}, where the objective function is plotted throughout the optimisations for the sharp as well as the rounded corner. We attribute this to the sensitivity filter and the continuous design variables, as this combination causes the area of intermediate design variables to decrease as the mesh is refined. Note that the coarse optimisation with a rounded corner seems to become infeasible and never recover, which might be attributed to numerical inconsistencies related to the sensitivity filter. 

\begin{figure*}[!htb]
\centering
\includegraphics[height=0.9\textheight]{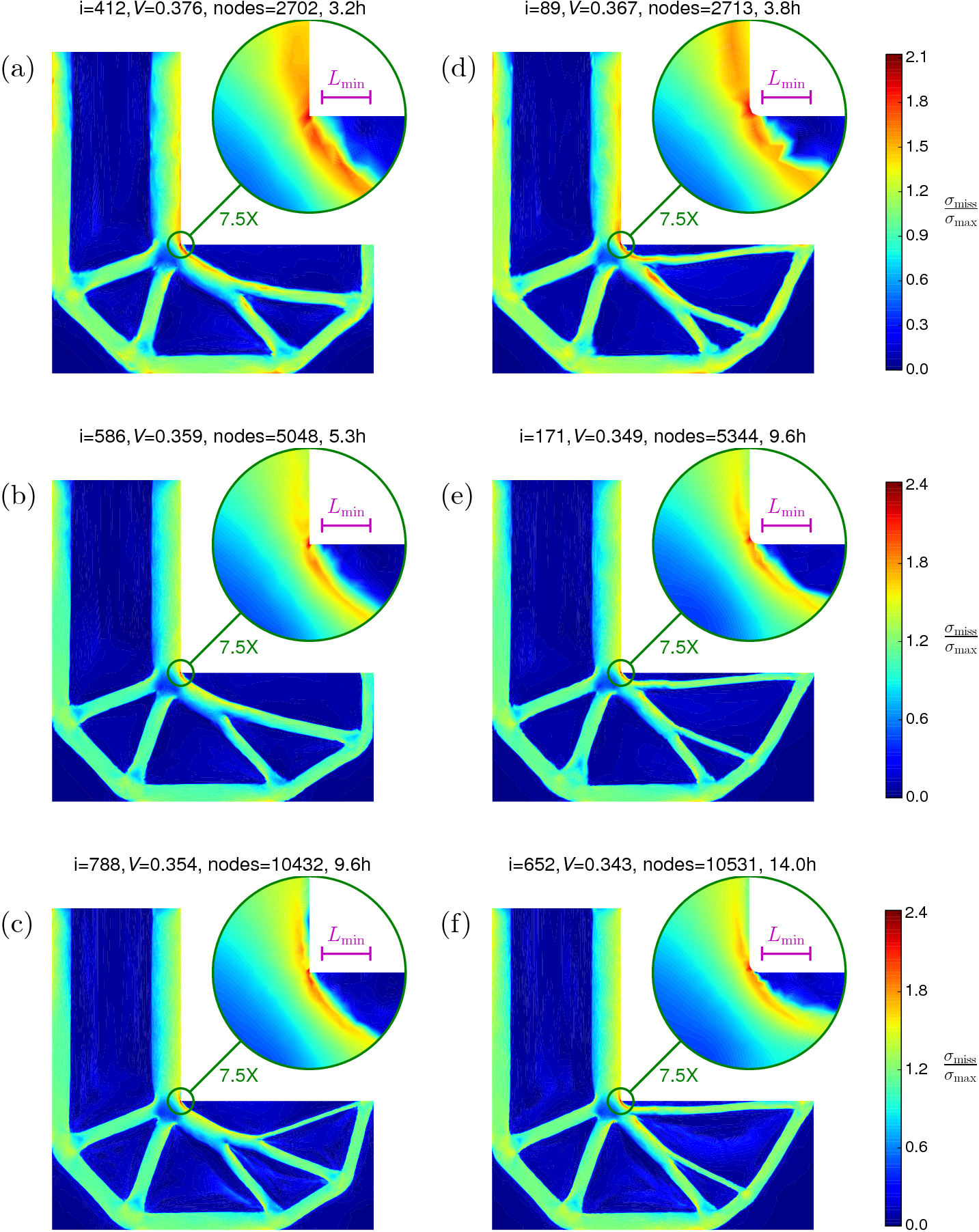}
\caption{Optimisations with a sharp (left) and a rounded corner (right) for $\eta_{\tilde{\rho}}$ equal to 0.03 (top), 0.015 (middle) and 0.0075 (bottom). The stress is shown for the iteration (i) at which the lowest volume fraction ($V$) occurs, while the stress and compliance constraints are satisfied. The stress calculated on the nodes, and the maximum values, in units of the maximum stress, are 2.11, 2.41 and 2.43 for the sharp corner and 2.05, 2.40 and 2.30 for the rounded corner. See online version for colours.}\label{fig:stress}
\end{figure*}

The stresses for the design in figure \ref{fig:designs} are plotted in figure \ref{fig:stress}. The maximum stress is not well resolved for $\eta_{\tilde{\rho}}=0.03$, but it does not change much between $\eta_{\tilde{\rho}}=0.015$ and $\eta_{\tilde{\rho}}=0.0075$. It seems like a safety margin no smaller than 2.5 is required, and as such figure \ref{fig:stress} highlights the need for post processing and verification of designs obtained using stress constrained topology optimisation.

Finally, the results of a preliminary study with an iterative solver is shown in table \ref{tab:tab}. The study indicates that the large element angles do not prohibit the use of iterative solvers, which is important for parallelism and efficiency in three dimensions.

\begin{table}[!htb]
\centering
\begin{tabular}{c | c c c c}
solver \textbackslash \,  $\eta_{\tilde{\rho}}$ & 0.03 & 0.015 & 0.0075 & NA\\ \hline
Direct (LU) & 0.4\,s & 0.72\,s & 1.47 & 0.49 \,s\\
Iterative (CG+AMG) & 0.61\,s & 1.27\,s & 2.69 & 0.68 \,s
\end{tabular}
\caption{The computational time for solving the forward problem using a direct solver (LU) for the designs with a rounded corner in figure \ref{fig:designs}. The computational time for solving the same linear systems with the conjugate gradient method and algebraic multigrid as preconditioner is also tabulated. We also show the same timings for the static mesh ($\eta_{\tilde{\rho}}=\mathrm{NA}$) in figure \ref{fig:static}. This demonstrates that the large element angles do not cause the relative performance of iterative solvers to decrease noticeably.}\label{tab:tab}
\end{table}

\paragraph*{Computational cost}
All computations are single threaded, but it is possible to perform mesh adaptation in parallel, see \cite{rokos2013thread, farrell2009galerkin}.

We terminate the optimisations using an iteration limit only, and the total computation time\footnote{Using an Intel(R) Core(TM) i7 870 @ 2.93GHz} before this triggers is shown in hours at the top of each plot in figure \ref{fig:designs}. We use a direct solver for the filter operations as well as for the forward and adjoint problems. The coarse optimisation of the L-bracket with a rounded corner probably only took around an hour to find the best design, but others have found similar designs with an order of magnitude less effort \cite{amstutz2010topological,bruggi2012topology} and comparing figures \ref{fig:static} and \ref{fig:designs}(a-c) reveals that it is possible to achieve a similar volume fraction using a fixed mesh and a similar amount of computational time. There is thus reason to believe that the presented method needs improvement, if it is to be competitive. It is, however, worth noting that the resolution of the design near the corner is significantly worse for the static mesh, even compared to the most coarse adaptive optimisation, so the fixed mesh might be benefiting from under resolving the stress singularity.

The optimisations with sharp corners are performed with PRAgMaTIc as mesh adaptation library, which is an optimised C++ implementation, so the actual mesh adaptation only takes 1-2\% of the total computation time. The optimisations with a rounded corner use an almost identical Octave/MATLAB implementation, which (although fully vectorised) is significantly slower, and therefore the mesh generation takes 20-30\% of the total computation time. Note that although the meshes conform well to the designs in figure \ref{fig:designs}, this is not true for the intermediate iterations. This is due to the fact that we optimise with $p=10$, so small imperfections in the designs can cause concentrations of stress away from the corner, which lead to many small elements in areas that are close to becoming void.  %In fact, optimisations with $\eta_{\tilde{\rho}}= 0.005$ often fail as the direct solvers runs out of memory. This behaviour is specific to stress constrained topology optimisation and could most likely be addressed by adaptively increasing the scaling factors for the metrics to impose a maximum element number.

We calculate the Steiner ellipse for each element, and use this to calculate the element aspect ratio by diving the radii product with the square of the smaller radius. This element aspect ratio is between 4 and 5 on average throughout the optimisations, which is a indication of speedup due to the anisotropic elements. One can expect the square of this in 3D corresponding to an order of magnitude reduction in computational cost compared to isotropic adaptation. 

\section{Conclusion}
We have performed stress constrained topology optimisations using a combination of anisotropic mesh adaptation and the method of moving asymptotes with interpolation of the asymptotes between iterations. We find that it is necessary to use continuous linear design variables and a sensitivity filter. The computational cost of introducing the mesh adaptation is negligible for an optimised implementation.

We argue that it might be beneficial to relax the problem statement to a rounded corner, as a radius of just 1\% of the characteristic length scale helps the optimiser find better designs. At least, in this case we are able to demonstrate mesh independence.

\section{Suggestions for Future Work}
It would be interesting to see the method applied to design of compliant mechanisms or meta materials due to the fact that the stress concentration results from the nature of the problem rather than the geometry.

Another interesting aspect is the post-processing of the result, i.e. to which extend is it beneficial to have many small elements aligned with the design contour, when the design is to be extracted?

Furthermore, the fact that anisotropic elements can give rise to large angles which cause problems for iterative solvers should be quantified in three dimensions and possibly also addressed by means of a meshing technique relying on advancing fronts \cite{loseille2014metric}. Another point for three dimensions is, that due to the trend of increasing number of cores in workstations, the fastest methods generally rely on shared memory parallelism \cite{christiansen2015combined}. Therefore one ought to investigate the benefit of parallelism, when extending the method to three dimensions.

Finally, the meshes that occur during the optimisation loop have only small variations during the later stages of the procedure, meaning that the degrees of freedom are chosen efficiently in the spatial dimension only, not in the ''optimisation dimension''. We suggest to address this issue using time-space elements and an optimiser defined at the continuous level.

\section{Acknowledgement}
This work is supported by the Villum Foundation.

%\bibliography{_Main}% Produces the bibliography via BibTeX.

\begin{thebibliography}{10}

\bibitem{habashi2000anisotropic}
Wagdi~G Habashi, Julien Dompierre, Yves Bourgault, Djaffar Ait-Ali-Yahia,
  Michel Fortin, and Marie-Gabrielle Vallet.
\newblock Anisotropic mesh adaptation: Towards user-independent,
  mesh-independent and solver-independent cfd. part i: General principles.
\newblock {\em International Journal for Numerical Methods in Fluids},
  32(6):725--744, 2000.

\bibitem{pain2001tetrahedral}
CC~Pain, AP~Umpleby, CRE De~Oliveira, and AJH Goddard.
\newblock Tetrahedral mesh optimisation and adaptivity for steady-state and
  transient finite element calculations.
\newblock {\em Computer Methods in Applied Mechanics and Engineering},
  190(29):3771--3796, 2001.

\bibitem{loseille2011continuous}
Adrien Loseille and Fr{\'e}d{\'e}ric Alauzet.
\newblock Continuous mesh framework part i: well-posed continuous interpolation
  error.
\newblock {\em SIAM Journal on Numerical Analysis}, 49(1):38--60, 2011.

\bibitem{loseille2010fully}
Adrien Loseille, Alain Dervieux, and Fr{\'e}d{\'e}ric Alauzet.
\newblock Fully anisotropic goal-oriented mesh adaptation for 3d steady euler
  equations.
\newblock {\em Journal of computational physics}, 229(8):2866--2897, 2010.

\bibitem{wallin2012optimal}
Mathias Wallin, Matti Ristinmaa, and Henrik Askfelt.
\newblock Optimal topologies derived from a phase-field method.
\newblock {\em Structural and Multidisciplinary Optimization}, 45(2):171--183,
  2012.

\bibitem{amstutz2010topological}
Samuel Amstutz and Antonio~A Novotny.
\newblock Topological optimization of structures subject to von mises stress
  constraints.
\newblock {\em Structural and Multidisciplinary Optimization}, 41(3):407--420,
  2010.

\bibitem{christiansen2015combined}
Asger~Nyman Christiansen, J~Andreas B{\ae}rentzen, Morten Nobel-J{\o}rgensen,
  Niels Aage, and Ole Sigmund.
\newblock Combined shape and topology optimization of 3d structures.
\newblock {\em Computers \& Graphics}, 46:25--35, 2015.

\bibitem{borrvall2001large}
Thomas Borrvall and Joakim Petersson.
\newblock Large-scale topology optimization in 3d using parallel computing.
\newblock {\em Computer methods in applied mechanics and engineering},
  190(46):6201--6229, 2001.

\bibitem{aagetopology}
Niels Aage, Erik Andreassen, and Boyan~Stefanov Lazarov.
\newblock Topology optimization using petsc.
\newblock {\em Structural and Multidisciplinary Optimization}, 2014.

\bibitem{aage2008topology}
Niels Aage, Thomas~H Poulsen, Allan Gersborg-Hansen, and Ole Sigmund.
\newblock Topology optimization of large scale stokes flow problems.
\newblock {\em Structural and Multidisciplinary Optimization}, 35(2):175--180,
  2008.

\bibitem{xia2012level}
Qi~Xia, Tielin Shi, Shiyuan Liu, and Michael~Yu Wang.
\newblock A level set solution to the stress-based structural shape and
  topology optimization.
\newblock {\em Computers \& Structures}, 90:55--64, 2012.

\bibitem{guo2011stress}
Xu~Guo, Wei~Sheng Zhang, Michael~Yu Wang, and Peng Wei.
\newblock Stress-related topology optimization via level set approach.
\newblock {\em Computer Methods in Applied Mechanics and Engineering},
  200(47):3439--3452, 2011.

\bibitem{zhang2013optimal}
Wei~Sheng Zhang, Xu~Guo, Michael~Yu Wang, and Peng Wei.
\newblock Optimal topology design of continuum structures with stress
  concentration alleviation via level set method.
\newblock {\em International Journal for Numerical Methods in Engineering},
  93(9):942--959, 2013.

\bibitem{bruggi2012topology}
Matteo Bruggi and Pierre Duysinx.
\newblock Topology optimization for minimum weight with compliance and stress
  constraints.
\newblock {\em Structural and Multidisciplinary Optimization}, 46(3):369--384,
  2012.

\bibitem{duysinx1998new}
Pierre Duysinx and Ole Sigmund.
\newblock New developments in handling stress constraints in optimal material
  distribution.
\newblock In {\em Proc of the 7th AIAA/USAF/NASAISSMO Symp on Multidisciplinary
  Analysis and Optimization}, volume~1, pages 1501--1509, 1998.

\bibitem{holmberg2013stress}
Erik Holmberg, Bo~Torstenfelt, and Anders Klarbring.
\newblock Stress constrained topology optimization.
\newblock {\em Structural and Multidisciplinary Optimization}, 48(1):33--47,
  2013.

\bibitem{farrell2013automated}
Patrick~E Farrell, David~A Ham, Simon~W Funke, and Marie~E Rognes.
\newblock Automated derivation of the adjoint of high-level transient finite
  element programs.
\newblock {\em SIAM Journal on Scientific Computing}, 35(4):C369--C393, 2013.

\bibitem{Note1}
Preliminary results related to this work was presented at FEniCS 14 in Paris,
  June 2014 and at the International Meshing Roundtable 23 in London, October
  2014.

\bibitem{loseille2014metric}
Adrien Loseille.
\newblock Metric-orthogonal anisotropic mesh generation.
\newblock {\em Procedia Engineering}, 82:403--415, 2014.

\bibitem{chen2007optimal}
Long Chen, Pengtao Sun, and Jinchao Xu.
\newblock Optimal anisotropic meshes for minimizing interpolation errors in
  𝐿\^{}$\{$𝑝$\}$-norm.
\newblock {\em Mathematics of Computation}, 76(257):179--204, 2007.

\bibitem{vasilevski2005error}
Yu~V Vasilevski and KN~Lipnikov.
\newblock Error bounds for controllable adaptive algorithms based on a hessian
  recovery.
\newblock {\em Computational Mathematics and Mathematical Physics},
  45(8):1374--1384, 2005.

\bibitem{rokos2013thread}
Georgios Rokos, Gerard~J Gorman, James Southern, and Paul~HJ Kelly.
\newblock A thread-parallel algorithm for anisotropic mesh adaptation.
\newblock {\em arXiv preprint arXiv:1308.2480}, 2013.

\bibitem{bendsoe2003topology}
Martin~Philip Bendsoe and Ole Sigmund.
\newblock {\em Topology optimization: theory, methods and applications}.
\newblock Springer, 2003.

\bibitem{cheng1997varepsilon}
GD~Cheng and Xiao Guo.
\newblock $\varepsilon$-relaxed approach in structural topology optimization.
\newblock {\em Structural Optimization}, 13(4):258--266, 1997.

\bibitem{svanberg1987method}
Krister Svanberg.
\newblock The method of moving asymptotes—a new method for structural
  optimization.
\newblock {\em International journal for numerical methods in engineering},
  24(2):359--373, 1987.

\bibitem{suresh2013stress}
Krishnan Suresh and Meisam Takalloozadeh.
\newblock Stress-constrained topology optimization: a topological level-set
  approach.
\newblock {\em Structural and Multidisciplinary Optimization}, 48(2):295--309,
  2013.

\bibitem{wang2011projection}
Fengwen Wang, Boyan~Stefanov Lazarov, and Ole Sigmund.
\newblock On projection methods, convergence and robust formulations in
  topology optimization.
\newblock {\em Structural and Multidisciplinary Optimization}, 43(6):767--784,
  2011.

\bibitem{le2010stress}
Chau Le, Julian Norato, Tyler Bruns, Christopher Ha, and Daniel Tortorelli.
\newblock Stress-based topology optimization for continua.
\newblock {\em Structural and Multidisciplinary Optimization}, 41(4):605--620,
  2010.

\bibitem{lazarov2011filters}
Boyan~Stefanov Lazarov and Ole Sigmund.
\newblock Filters in topology optimization based on helmholtz-type differential
  equations.
\newblock {\em International Journal for Numerical Methods in Engineering},
  86(6):765--781, 2011.

\bibitem{LoggMardalEtAl2012a}
Anders Logg, Kent-Andre Mardal, Garth~N. Wells, et~al.
\newblock {\em Automated Solution of Differential Equations by the Finite
  Element Method}.
\newblock Springer, 2012.

\bibitem{sigmund2012sensitivity}
Ole Sigmund and Kurt Maute.
\newblock Sensitivity filtering from a continuum mechanics perspective.
\newblock {\em Structural and Multidisciplinary Optimization}, 46(4):471--475,
  2012.

\bibitem{farrell2009galerkin}
Patrick~E Farrell.
\newblock {\em Galerkin projection of discrete fields via supermesh
  construction}.
\newblock PhD thesis, Imperial College London, 2009.

\bibitem{Note2}
Using an Intel(R) Core(TM) i7 870 @ 2.93GHz.

\end{thebibliography}
\providecommand{\noopsort}[1]{}\providecommand{\singleletter}[1]{#1}%

\end{document}